\documentclass[11pt, letterpaper, oneside]{amsart}

\usepackage{latexsym, amsmath, amssymb, amsfonts, amscd}
\usepackage{amsthm}
\usepackage{t1enc}
\usepackage[mathscr]{eucal}
\usepackage{indentfirst}
\usepackage{graphicx, pb-diagram}
\usepackage{fancyhdr}
\usepackage{fancybox}
\usepackage{enumerate}
\usepackage[all]{xy}	 	
\usepackage{url}

\theoremstyle{plain}
\newtheorem{thm}{Theorem}[section]

\newtheorem{prop}[thm]{Proposition}
\newtheorem{lemma}[thm]{Lemma}
\newtheorem{cor}[thm]{Corollary}

\theoremstyle{definition}
\newtheorem{defi}[thm]{Definition}

\theoremstyle{remark}
\newtheorem{remark}[thm]{Remark}
\newtheorem{ex}[thm]{Example}

\newcommand{\CC}{\ensuremath{\mathbb C}}
\newcommand{\RR}{\ensuremath{\mathbb R}}

\newcommand{\g}{\ensuremath{\frak{g}}}

\newcommand{\un}{\underline}

\newcommand{\cL}{\mathcal{L}}
\newcommand{\cG}{\mathcal{G}}

\newcommand{\cF}{\mathcal{F}}
\newcommand{\cD}{\mathcal{D}}
\newcommand{\cJ}{\mathcal{J}}
\newcommand{\kp}{K^{\perp}}

\newcommand{\he}{\hat{e}}

\newcommand{\te}{\tilde{e}}

\newcommand{\Gbw}{\Gamma_{bas}(\kp)}
\newcommand{\ra}{\rangle}
\newcommand{\la}{\langle}
\newcommand{\Gbj}{\Gamma_{bas}(\kp\cap \cJ \kp )}
\newcommand{\kpj}{\kp\cap \cJ \kp}

\begin{document}

\title{Reduction of branes in generalized complex geometry}
\author{  Marco Zambon}
\address{Institut f\"ur Mathematik, Universit\"at Z\"urich-Irchel, Winterthurerstr. 190, CH-8057 Z\"urich, Switzerland}
\email{marco.zambon@math.unizh.ch}

\begin{abstract}
We show that  certain submanifolds of  generalized complex manifolds
(``weak branes'') admit a natural
quotient which inherits a generalized complex structure. 
This is analog to quotienting coisotropic submanifolds of symplectic manifolds.
In particular Gualtieri's generalized complex submanifolds (``branes'')
quotient to space-filling branes.
Along the way we perform reductions by foliations (i.e. no group action is involved)
for exact Courant algebroids - interpreting the reduced \v{S}evera class -
and for Dirac structures.
\end{abstract}

\maketitle
\tableofcontents

\section{Introduction}\label{intro}
Consider the following setup in ordinary geometry: a manifold $M$ and a submanifold $C$ endowed with some
integrable distribution $\cF$ so that $\un{C}:=C/\cF$ be smooth. Then we have a projection
$pr:C \rightarrow \un{C}$ which induces a vector bundle morphism 
$pr_*:TC \rightarrow T\un{C}$. If $M$ is endowed with some geometric structure,
such as a symplectic 2-form $\omega$,  one can ask when $\omega$ induces
 a symplectic form on $\un{C}$. 

This happens for example when $C$ is a
coisotropic submanifold\footnote{This means that the symplectic orthogonal
of $TC$ is contained in $TC$.}. Indeed in this case the pullback $i^*\omega$ 
of $\omega$ to $C$ has a kernel $\cF$ which is of constant rank and integrable, and the closeness of
$\omega$ ensures that if $p$ and $q$ lie in the same $\cF$-leaf
then
$(i^*\omega)_p$ 
and $(i^*\omega)_q$ project  to the same linear symplectic form at $pr(p)=pr(q)$, so that one obtains
a well-defined
symplectic form on $\un{C}$.
An instance of the above is when there is a Lie group $G$ acting 
hamiltonianly on $M$ with moment map $\nu:M \rightarrow  {\g}^*$ and  $C$ is the zero
level set of $\nu$ (Marsden-Weinstein reduction \cite{MW}).\\

In this paper we consider the geometry that arises when one replaces the 
tangent bundle $TM$ with an exact Courant algebroid $E$ over $M$ (any such $E$ is non-canonically
isomorphic to $TM \oplus T^*M$).
In this context reduction by the action of a Lie group
has been considered by several
authors (Bursztyn-Cavalcanti-Gualtieri \cite{BCG},  Hu \cite{Hu1,Hu2},
Stienon-Xu \cite{SX},
Tolman-Lin \cite{LT2,LT1}); in this paper we do not assume any group action.
Unlike the tangent
bundle case, knowing  $\un{C}$ does not automatically determine the exact Courant algebroid over it.
We have to replace
 the foliation $\cF$ by more data, namely a 
  suitable subbundle $K$ of $E|_C$, and we construct by a quotienting procedure a Courant algebroid $\un{E}$ on $\un{C}$ 
 (Theorem \ref{redcour}).  Our construction follows closely the one of
 \cite{BCG}, where
the  group action provides an identification between  fibers of $E$ at different points;
 in our case we make up for this by asking that there exist enough ``basic sections'' (Def. \ref{bas}).
 Further, we describe in a simple way (see Def. \ref{adapted})
 which splittings of $E$ induce 3-forms on $M$  (representing 
 the \v{S}evera class of $E$) which descend to  3-forms on $\un{C}$   (representing 
 the \v{S}evera class of $\un{E}$).
 
Once we know how to reduce  an exact Courant algebroid, we can ask when geometric structures defined on them
descend to the reduced  exact Courant algebroid. We consider Dirac structures
(suitable subbundles of $E$) and generalized complex structures (suitable endomorphisms of $E$).
We give sufficient conditions for these
structures to descend in Prop. \ref{reddirac} and Prop. \ref{redgcs} respectively.
The  ideas and techniques are borrowed the literature cited above, in particular from 
\cite{BCG} and \cite{SX} (however our proof differs from these two references in that
 we reduce generalized complex structures
directly and not viewing them as Dirac structures in the complexification of $E$).

The heart of this paper is Section \ref{branes}, where we identify the objects that automatically
satisfy the assumptions needed to perform generalized complex reduction. When $M$ is a generalized
complex manifold we consider pairs consisting of a submanifold $C$ of $M$ and suitable maximal
isotropic subbundle $L$ of $E|_C$ (we call them ``weak branes'' in Def. \ref{defweak}).
We show in Prop. \ref{weakbrane} that weak branes admit a canonical quotient $\un{C}$ which is endowed with an exact Courant algebroid
and a generalized complex structure; this construction is inspired by Thm. 2.1 of Vaisman's work
\cite{Vared} in the setting of the standard Courant algebroid.

Particular cases of weak branes are generalized complex submanifolds $(C,L)$ (also known as ``branes'',
see Def. \ref{gcsub}), which were first introduced by Gualtieri \cite{Gu} and are relevant to physics
 \cite{KO}. Using our reduction of  Dirac structures
we show in Thm. \ref{redbrane} that the quotients $\un{C}$ of branes, which by the above are generalized complex
manifolds, are also endowed with the  structure of a space-filling brane.   This is interesting also because space filling branes
induce an honest complex structure on the underlying manifold \cite{Gu2}.\\

The reduction statements we had to develop 
in order to prove the results of  Section \ref{branes}  are
versions ``without group action'' of statements that
already appeared in the literature \cite{BCG}\cite{Hu1,Hu2}
 \cite{ BS,SX}
  \cite{LT2,LT1} \cite{Vared}. Consequently many ideas and techniques are borrowed from the existing
  literature; we make appropriate references in the text whenever possible.
  In particular we followed closely \cite{BCG} (also as far as notation and conventions are
  concerned). \\
 
\textbf{Plan of the paper:} in Section \ref{review} we review exact Courant algebroids. In Section
 \ref{cour} we perform the reduction of exact Courant algebroids and
determine objects that naturally satisfy the assumptions needed for the reduction. In Section \ref{dirac} we 
perform the reduction of Dirac structures.
In section \ref{Sev} we describe  the reduced \v{S}evera class. 
In Section \ref{gcs} we reduce generalized complex structures and comment briefly on generalized K\"ahler
reduction. The main section of this paper is Section \ref{branes}: we reduce branes and weak branes,
providing few examples. We also give a criteria that allows to obtain weak branes by restricting
to cosymplectic submanifolds.\\
 
\textbf{Acknowledgments:} I am very indebted to Henrique Bursztyn for clarifying to me
some of the constructions of \cite{BCG}, and to Marco Gualtieri for 
some crucial discussions at the Geometry Conference in Honor of Nigel Hitchin (September 2006), during which
he clarified conceptual issues and also suggested some of the techniques used in this paper. 
Furthermore I thank A. Cattaneo, G. Cavalcanti, F. Falceto, B. Scardua and M. Stienon
for discussions. I also thank A. Cattaneo for supporting my attendance to   conferences 
relevant to this work
and
 F. Falceto for a visit to the Universidad de Zaragoza where part of this work was
done.

 Further I acknowledge support from the Forschungskredit of the Universit\"at Z\"urich
and partial support of SNF Grant No.~20-113439.
This work has been partially supported
by the European Union through the FP6 Marie Curie RTN ENIGMA (Contract
number MRTN-CT-2004-5652) and by the European Science Foundation
through the MISGAM program. 

\section{Review of Courant algebroids}\label{review}

We review the notion of exact Courant algebroid; see \cite{BCG} and \cite{Hu1} for more details.

\begin{defi}\label{ca}
A {\it Courant algebroid} over a manifold $M$ is a vector bundle $E
\to M$ equipped with a fibrewise non-degenerate symmetric bilinear
form $\la \cdot,\cdot \ra$, a bilinear bracket $[\cdot,\cdot]$ on
the smooth sections $\Gamma(E)$, and a bundle map $\pi: E\to TM$
called the \textit{anchor}, which satisfy the following conditions
for all $e_1,e_2,e_3\in \Gamma(E)$ and $f\in C^{\infty}(M)$:
\begin{itemize}
\item[C1)] $[e_1,[e_2,e_3]] = [[e_1,e_2],e_3] +
[e_2,[e_1,e_3]]$,
\item[C2)] $\pi([e_1,e_2])=[\pi(e_1),\pi(e_2)]$,
\item[C3)] $[e_1,fe_2]=f[e_1,e_2]+ (\pi(e_1)f) e_2$,
\item[C4)] $\pi(e_1)\langle e_2,e_3 \rangle= \langle [e_1,e_2],e_3 \rangle
+ \langle e_2, [e_1, e_3]\rangle$,
\item[C5)] $[e_1,e_1] = \cD \langle e_1,e_1 \rangle$,
\end{itemize}
where $\cD =\frac{1}{2}\pi^*\circ d: C^{\infty}(M)\rightarrow
\Gamma(E)$ (using $\langle \cdot,\cdot \rangle$ to identify $E$
with $E^*$).
\end{defi}

We see from axiom $C5)$ that the bracket is not skew-symmetric:
$$
[e_1,e_2]=-[e_2,e_1]+2\cD \langle e_1,e_2 \rangle.$$ Hence we have
the following ``Leibniz rule for the first entry'':
$[fe_1,e_2]=f[e_1,e_2]- (\pi(e_2)f) e_1+2\langle e_1,e_2 \rangle \cD
f$.

\begin{defi} A Courant algebroid is {\it exact} if
the following sequence is exact:
\begin{equation}\label{exact}
0 \longrightarrow T^*M \stackrel{\pi^*}{\longrightarrow} E
\stackrel{\pi}{\longrightarrow} TM \longrightarrow 0
\end{equation}
\end{defi}
 To simplify the notation, 
in the sequel we will often omit the map $T^*M
\overset{\pi^*}\rightarrow E^* \cong E$ and think of $T^*M$ as being
a subbundle of $E$.
Given an exact Courant algebroid, we may always choose a right
splitting $\sigma: TM\to E$   whose image in $E$ is isotropic with
respect to $\langle \cdot,\cdot \rangle$. Such a splitting induces
the closed 3-form on $M$  given by 
$$H(X,Y,Z)=2 \la [\sigma X,\sigma
Y], \sigma Z \ra.$$
 Using
the bundle isomorphism $\nabla+\tfrac{1}{2}\pi^*:TM\oplus T^*M \to
E$, one can transport the Courant algebroid structure onto $TM\oplus
T^*M$.  The resulting structure is as follows 
(where $X_i+\xi_i \in \Gamma(TM\oplus TM^*)$): the bilinear pairing is
\begin{equation}\label{pairing}
\langle X_1+\xi_1, X_2+\xi_2 \rangle = \frac{1}{2}(\xi_2(X_1) + \xi_1(X_2)),
\end{equation}
and the bracket is
\begin{equation}\label{Hcour}
[ X_1+\xi_1, X_2+\xi_2]_H= [X_1,X_2] + \cL_{X_1} \xi_2 - i_{X_2} d\xi_1  + i_{X_2} i_{X_1} H,
\end{equation}
which is the $H$-twisted Courant bracket on $TM\oplus T^*M$
\cite{SW}.  Isotropic splittings of $\eqref{exact}$ differ by
2-forms $b\in \Omega^2(M)$, and a change of splitting modifies the
curvature $H$ by the exact form $db$. Hence there is a well-defined 
 cohomology class
$[H]\in H^3(M,\mathbb{R})$
attached to the exact Courant
algebroid structure on $E$; $[H]$ is called the \textit{\v{S}evera class} of $E$.
 
We refer to  \cite{BCG} and \cite{Hu1} for  information on 
the group of automorphisms $Aut(E)$
and its Lie algebra $Der(E)$. Here we just mention few facts, the
first of which underlies many of our constructions:
 for any 
$e\in \Gamma(E)$, 
  $[e,\cdot]$ is an element of $Der(E)$ and hence integrates to an automorphism of the
  Courant algebroid $E$.  
 Notice that for closed 1-forms $\xi$ (seen as sections of 
  $T^*M\subset E$)
  we have $[\xi,\cdot]=0$ by \eqref{Hcour}.    Further, any
   2-form $B$ on $M$ determines a vector bundle map $TM\oplus TM^*\rightarrow TM\oplus TM^*$ by 
$e^B:X+\xi\mapsto X+\xi+i_XB$ \cite{Gu} and these ``gauge
transformations'' satisfy
\begin{equation}\label{gauge}
[e^B\ \cdot, e^B\ \cdot]_H = e^B[\cdot,\cdot]_{H+dB}.
\end{equation}

\section{The case of exact Courant algebroids}\label{cour}

In this section we reduce exact Courant algebroids (Thm. \ref{redcour}).
  
Let $M$ be a manifold, $E$ an exact Courant algebroid over $M$, and $C$ a
 submanifold.

\begin{lemma}\label{bra}
Let $D\rightarrow C$ be a subbundle of $E$ such that
$\pi(D^{\perp})\subset TC$ (where $D^{\perp}$ denotes the orthogonal
to $D$ w.r.t. the symmetric pairing), and $e_1,e_2$ sections of
$D^{\perp}$. Then the expression $[\te_1,\te_2]|_C$, where $\te_i$
are extensions of $e_i$ to sections of $E\rightarrow M$, depends on
the extensions only up to sections of $D$.
\end{lemma}
\begin{proof}
Fix extensions $\te_i$ of $e_i$ ($i=1,2$). We have to show that for
functions $f_i$ vanishing on $C$ and sections $\he_i$ of $E$ we have
$[\te_1+f_1 \he_1,\te_2+f_2 \he_2]|_C=[\te_1,\te_2]|_C$ up to
sections of $D$. By the Leibniz rule $C3)$ 
and since $\pi(e_1)\subset TC$
we have
$[\te_1,f_2 \he_2]|_C=0$. Also $[f_1 \he_1,\te_2]|_C=2\langle
\he_1,\te_2 \rangle (\cD f_1)|_C\subset
N^*C\subset (\pi(D^{\perp}))^{\circ}=D\cap T^*M$
\footnote{Indeed for any subspace $D$ of a
vector space $T\oplus T^*$,  we have $D\cap T^*=(\pi (D^{\perp}))^{\circ}$.}.
 The term $[f_1 \he_1,f_2
\he_2]|_C$ vanishes by the above since $(f_1 \he_1)|_C$ is a section
of $D$.
\end{proof}

\begin{remark}\label{D}
If $D\rightarrow C$ is a subbundle of $E$ such that
$\pi(D^{\perp})\subset TC$ we can make sense of a statement like
``$[e_1,e_2]\subset D$'' for $e_1,e_2\in \Gamma(D^{\perp})$: it means that
$[\te_1,\te_2]|_C\subset D$ for one (or equivalently, by
Lemma \ref{bra}, for all) extensions $\te_i$ to sections of
$E\rightarrow M$. Similarly, we take $[\Gamma(D^{\perp}),\Gamma(D^{\perp})]\subset
\Gamma(D)$ to mean $[e_1,e_2]\subset D$ for all
$e_1,e_2\in  \Gamma(D^{\perp})$.
\end{remark}

Now fix an \emph{isotropic} subbundle $K\rightarrow C$ of $E$, i.e.
$K\subset \kp$, such that $\pi(\kp)=TC$.

\begin{defi}\label{bas}
We define the space of sections of $\kp$ which are \emph{basic w.r.t. $K$}
as
\begin{equation}
\Gamma_{bas}(\kp):=\{e\in \Gamma(\kp):[\Gamma(K),e]\subset
\Gamma(K)\}.
\end{equation}
\end{defi}

\begin{remark}\label{localk}
To ensure that a section  $e$ of $\kp$ be basic it suffices to consider locally defined sections
of $K$ that span $K$ point-wise. That is, it suffices to 
show that
for every point of $C$ there is a neighborhood $U\subset C$ and a
 subset $S\subset \Gamma(K|_U)$ with 
  $span\{k_p:k\in S\}=K_p$ (for every $p\in U$)
 so that
$[S,e|_U]\subset \Gamma(K|_U)$.
Indeed from the ``Leibniz rule in the first entry'' it follows that
$[\Gamma(K),e]\subset \Gamma(K)$ .
\end{remark}

\begin{lemma}\label{subset}
Assume that  the sections
of $\Gbw$ span $\kp$ at every point, i.e. that $span\{e_p:e\in \Gbw\}=\kp_p$ for every $p\in C$. Then
\itemize
\item[1)] $[\Gamma(K),\Gamma(\kp)]\subset \Gamma(\kp)$
\item[2)] $[\Gamma(K),\Gamma(K)]\subset \Gamma(K)$.
\end{lemma}
\begin{proof}
Fix a subset of sections   
$\{e_i\} \subset \Gbw$ that spans point-wise $\kp$. For any section $k$ of $K$ and
functions $f_i$ (so that the sum $\sum f_i
e_i$ is locally finite) by the Leibniz rule  we have $[k, \sum f_i
e_i]\subset \kp$, proving 1). Now 1) is equivalent to 2), as
can be seen using axiom C4) in the definition of Courant algebroid:
let $k_1,k_2$ be   sections of $K$ and $e$ a section of $\kp$. Then
$\la[k_1,e],k_2 \ra+\la e,[k_1,k_2] \ra=\pi (k_1) \la e,k_2 \ra=0$
because $\pi (K)\subset \pi (\kp) =TC$.
\end{proof}

\begin{remark}
A converse to Lemma \ref{subset} for local sections is given in \cite{CFZ}.
\end{remark}

The proof of the following theorem is modeled on Thm. 3.3 of \cite{BCG}.
When referring to the smoothness of the quotient of a manifold by a foliation, we refer to the unique differentiable  structure so that the projection map is a submersion.  
\begin{thm}[Exact Courant algebroid reduction]\label{redcour}
 Let
$E$ be an exact Courant algebroid over $M$, $C$ a submanifold of
$M$, and $K$ an isotropic subbundle of $E$ over $C$ such that
$\pi(\kp)=TC$. Assume that the space of (global) sections 
 $\Gbw$ spans point-wise $\kp$
 (i.e. that  $span\{e_p:e\in \Gbw\}=\kp_p$ for every $p\in C$)
    and that the quotient $\un{C}$ of $C$ by the foliation
integrating $\pi(K)$ be a smooth manifold. Then there is an exact
Courant algebroid $\un{E}$ over $\un{C}$ that fits in the following 
pullback diagram of vector bundles:
\[\xymatrix{  \kp/K \ar[r]\ar[d]  & \un{E}\ar[d]
\\
C \ar[r]& \un{C}}.
\]
\end{thm}

\begin{proof} Notice that since
 $\pi (K)$ has  constant rank iff $\pi (D^{\perp})$ does (use the previous footnote
 or eq. (2.17) of \cite{Vaiso}) it follows that $\pi(K)$ is a  
 regular distribution on $C$. Further, by the assumption on basic sections
 and item 2) of Lemma \ref{subset},
$\pi (K)$ is an  integrable distribution, so there exists a regular foliation
integrating $\pi(K)$. We divide the proof in 2 steps.

\textbf{Step 1}
 To describe   the
vector bundle $\un{E}$ we have to explain how we identify fibers of
$\kp/K$ over two points $p,q$ lying in the same leaf $F$ of $\pi (K)$.
We do this as follows: we identify two elements  $\hat{e}(p)\in
(\kp/K)_p$ and $\hat{e}(q)\in(\kp/K)_q$ iff there is a section $e\in
\Gbw$ which under the projection $\kp\rightarrow \kp/K$ maps\footnote{In other words, we give a canonical
trivialization of $(\kp/K)|_F$ by projecting into it a frame for
$\kp|_F$ consisting of basic sections; by assumptions we have enough
basic sections to really get a frame for $(\kp/K)|_F$.}
 to  
$\hat{e}(p)$ at $p$ and $\hat{e}(q)$ at $q$.
 To show that this procedure gives a well-defined identification of
 $(\kp/K)_p$ and $(\kp/K)_q$, we need to show that if $e_1$ and $e_2$ are sections
 of $\Gbw$ such
that $e_1(p)$ and $e_2(p)$ map to $\hat{e}(p)$, then $e_1(q)$ and
$e_2(q)$ map to the same element of $(\kp/K)_q$.

 Pick a finite sequence of local sections
$k_1,\dots,k_n$ of $K$ that join $p$ to $q$, i.e. such that following
successively the vector fields $\pi (k_i)$ for times $t_i$ the point
$p$ is mapped to $q$. Extend each $k_i$ to a section $\tilde{k}_i$
of $E$. Denote by $e^{ad_{\tilde{k}_i}}$  the Courant algebroid
automorphism of $E$ obtained integrating
$ad_{\tilde{k}_i}=[\tilde{k}_i,\cdot]$, and by $\Phi$ the
composition $e^{ad_{t_n\tilde{k}_n}}\circ\dots \circ
e^{ad_{t_1\tilde{k}_1}}$. Since $e_1$ is a basic section 
we have
$[k_i,e_1]\subset K$ for all $i$. So
$\Phi (e_1(p))-e_1(q)\in K_q$, and similarly for $e_2$. Now
$e_1(p)-e_2(p)\in K_p$ by assumption, so
 because of item 2) of Lemma \ref{subset}
 we have $\Phi (e_1(p)-e_2(p))\in K_q$.
We deduce that $e_1(q)-e_2(q)$ also belong to $K_q$ and
therefore project to the zero vector in $(\kp/K)_q$.

It is clear that   $\un{E}$, obtained from ${\kp}/K$ by identifying
the fibers over each leaf of $\pi( K)$ as above, is endowed with a
projection
 $\un{pr}$ onto $\un{C}$ (induced from the
projection $pr:{\kp}/K\rightarrow C$). $\un{E}$ is indeed a smooth
vector bundle: given any point $\un{p}$ of $\un{C}$ choose a
preimage $p\in C$ and a submanifold $S\subset C$ through $p$ transverse to
the leaves of $\pi (K)$. $S$ provides a chart around $\un{p}$ for
the manifold  $\un{C}$, and $pr^{-1}(S)$ is a vector subbundle  of
$\kp/K$ proving a chart for $\un{E}$ around $\un{p}$.

Notice that pulling back by the vector bundle epimorphism $\kp/K \rightarrow
\un{E}$ we can embed the space of sections of $\un{E}$ into the
space of sections of $\kp/K$, the image being the image of $\Gbw$
under the map $\kp \rightarrow \kp/K$. In
other words, we have a canonical identification $\Gamma(\un{E})\cong
\Gbw/\Gamma(K)$.

\textbf{Step 2}
The pairing $\la \cdot,\cdot \ra$ on  the fibers of $E$ induces
 a non-degenerate symmetric bilinear form on each fiber of  $\kp/K$, which  moreover 
 descends to $\un{E}$, because for any two
given sections $e_1,e_2 \in \Gbw$  $k \in \Gamma(K)$ we
have $\pi (k) \la e_1,e_2 \ra=0$ using C4).

For the bracket of sections of $\un{E}$, first notice that $\Gbw$
is closed (in the sense of Remark \ref{D}) under the bracket $[\cdot,\cdot]$ of $E$: if $e_1,e_2 \in
\Gbw$,  $[e_1,e_2]$ is a
section of $\kp$ by the argument (using
  C4)) in  Thm. 3.3 of \cite{BCG}.
Further $[e_1,e_2]$  is again basic by the ``Jacobi identity'' C1): for any
section $k$ of $K$ we have
$[k,[e_1,e_2]]=[[k,e_1],e_2]+[e_1,[k,e_2]]$. Now by definition of
basic section each $[k,e_i]$ lies in $K$, and applying once more the
definition of basic section\footnote{Together with the fact that for
any section $\hat{k}$ of $K$ we have
$[e_1,\hat{k}]=-[\hat{k},e_1]+2\cD \la e_1,\hat{k} \ra$
and $\cD \la e_1,\hat{k} \ra\subset N^*C=K\cap T^*M$.} we see that
$[k,[e_1,e_2]]\subset K$, i.e. that $[e_1,e_2]$ is basic.
In the light of Lemma \ref{bra}, what we really have a well-defined
bilinear form $\Gbw \times\Gbw \rightarrow \Gbw/\Gamma(K)$. Using the definition of
basic section we then have an induced bracket on
$\Gbw/\Gamma(K)$, which as we saw is canonically isomorphic to
$\Gamma(\un{E})$.

The induced anchor map $\un{\pi}: \un{E}\rightarrow T\un{C}$ well-defined since
$\pi(e)$ is a projectable vector field for
any basic section $e \in \Gbw$, as follows using  C2).

It is straightforward to check that 
$\un{E}\rightarrow
\un{C}$, endowed with the induced   symmetric
pairing, bilinear bracket on $\Gamma(\un{E})$ and  anchor
$\un{\pi}$, satisfies 
  axioms C1)-C5)
in the definition of Courant algebroid (Def. \ref{ca}).
Further the proof of Thm. 3.3 of \cite{BCG} shows that $\un{E}$ is an exact Courant algebroid. 

\end{proof}


\begin{remark}\label{morph} 
  The subbundle $\{(e,\un{e})|e \in \kp\} \rightarrow   \{(p,\un{p})|p \in C\}$  of $(E\times 
\un{E})\rightarrow M\times \un{C}$  provides a morphism of Courant algebroids from $E$ to $\un{E}$
as defined (essentially) in Def. 6.12 of \cite{AX} or in Def. 3.5.1 of \cite{HU}.\end{remark}

\begin{ex}[Quotients of submanifolds]\label{Cf}
Take $E$ to be $TM\oplus T^*M$ with the untwisted bracket. Let $C$ be a submanifold
endowed with a regular  distribution $\cF$, and
assume that the quotient $\un{C}=C/\cF$ be smooth. Take
$K:=\cF\oplus N^*C$.  $\Gamma(K)$ is spanned by
vector fields on $C$ lying in $\cF$ and differentials of
functions vanishing on  $C$. Since the latter 
act trivially, it is enough to consider the action of a vector
field $X\subset \cF$. Let $Y\oplus df|_C$ be a section of $\kp$, where
$Y$ is a projectable vector field  and  $f$ is the extension to $M$
of the pullback of a function on $\un{C}$. The action of $X$ on this
section is just $[X,Y]\oplus (\cL_X df)|_C$, which lies again in $K$.
Since such $Y\oplus df|_C$ span $\kp$ we can apply Thm. \ref{redcour}
and obtain a reduced Courant algebroid on $\un{C}$, namely
  $T\un{C}\oplus T^*\un{C}$ with the untwisted bracket.
\end{ex}

\begin{ex}
Let $E$ be an exact Courant algebroid over $M$ and $C$ a
submanifold of $M$. Then with $\hat{K}=N^*C$   the assumptions of Thm. \ref{redcour}
are satisfied; indeed all the sections of $\hat{K}^{\perp}=\pi^{-1}(TC)$ are
basic. Hence we recover Lemma 3.7 of \cite{BCG}.
\end{ex}

\section{The case of Dirac structures}\label{dirac}

Let $E$ be an exact Courant algebroid over $M$. Recall \cite{Cou} that a
\emph{Dirac structure}  is a maximal isotropic  
subbundle of $E$ which is closed under the Courant bracket. Now we let $C$
be a submanifold of $M$ and consider a maximal isotropic
subbundle $L\subset E$ defined over $C$ (not necessarily satisfying
 $\pi(L)\subset TC$). The following is analog to Thm. 4.2 of \cite{BCG}.
\begin{prop}[Dirac reduction]\label{reddirac}
Let $E\rightarrow M$ and  $K\rightarrow C$ satisfy the assumptions
of Thm. \ref{redcour}, so that we have an exact Courant algebroid
$\un{E}\rightarrow \un{C}$. Let $L$ be a maximally isotropic
  subbundle of $E|_C$ such that $L\cap \kp$ has constant rank,
  and assume that  \begin{eqnarray}\label{diracdes}
[\Gamma(K),\Gamma(L\cap \kp)]\subset \Gamma(L+K).
 \end{eqnarray}
  Then
$L$ descends to a maximal isotropic subbundle $\un{L}$ of
$\un{E}\rightarrow \un{C}$. If furthermore \begin{eqnarray}\label{diracint}
[\Gamma_{bas}(L\cap \kp), \Gamma_{bas}(L\cap \kp)] \subset
\Gamma(L+K).
 \end{eqnarray}
 then $\un{L}$ is an
(integrable) Dirac structure. Here $\Gamma_{bas}(L\cap \kp):=\Gamma(L)\cap \Gbw$
\end{prop}
\begin{proof}
At every $p\in C$ we have a Lagrangian relation
between 
 $E_p$ and
$(\kp/K)_p$ given by $\{(e,e+K_p): e\in \kp_p\}$. The image of $L_p$ under this relation, which we denote
by $\un{L}(p)$, is  maximal isotropic because $L_p$ is. Doing
this at every point of $C$ we obtain a maximally isotropic subbundle
of $\kp/K$, which is furthermore smooth because  $L(p)$ is the image
of $(L\cap \kp)_p$, which has constant rank by assumption,  under
the projection $\kp_p\rightarrow (\kp/K)_p$.

Recall that in Thm. \ref{redcour} we identified $  (\kp/K)_p$ and
$(\kp/K)_q$ when $p$ and $q$ lie in the same leaf of $\pi(K)$,
and that the identification was induced by the Courant algebroid
automorphism $\Phi$ of $E$ obtained integrating
  any sequence of locally defined sections $k_1,\dots,k_n$ of $K$ that join
$p$ to $q$ (see Remark \ref{alter}). Assumption \eqref{diracdes}
(with Lemma \ref{subset} 1)) is exactly what is needed to
ensure that $\Phi$ maps $L\cap \kp$ into
$(L+K)\cap \kp=(L\cap \kp)+K$, so that
   $\un{L}(p)$ gets identified with  $\un{L}(q)$. As
a consequence we obtain a well-defined smooth maximally isotropic
subbundle $\un{L}$ of the reduced Courant algebroid $\un{E}$, i.e.
an almost Dirac structure for $\un{E}$. Now assume that
\eqref{diracint} holds, and take two sections of $\un{L}$, which by
abuse of notation we denote $\un{e}_1,\un{e}_2$. Since $L\cap \kp$
has constant rank we can lift them to   sections $e_1,e_2$   of
$\Gamma_{bas}(L\cap \kp)$. As for all elements of $\Gbw$ their bracket
lies in $\Gbw$, and by assumption it also lies in $L+K$, so
$[e_1,e_2]$ is a basic section of $(L+K)\cap \kp=(L\cap \kp)+K$. Its
projection under $\kp/K\rightarrow \un{E}$, which is by definition
the bracket of $\un{e}_1$ and $\un{e}_2$, lies  then in $\un{L}$.
\end{proof}

\begin{ex}[Coisotropic reduction] \label{coiso}
Let $(M,\Pi)$ be a Poisson manifold and $C$ a coisotropic
submanifold\footnote{This means that $\sharp N^*C\subset TC$, where
$\sharp:T^*M\rightarrow TM$ is the contraction with $\Pi$.}. It is known \cite{Cou} that
the characteristic distribution $\cF:=\sharp N^*C$ is a singular
integrable distribution; assume that it is regular and the quotient
$\un{C}=C/\cF$ be smooth. It is known that $D=\{(\sharp \xi,\xi):\xi \in
T^*P\}$ is a Dirac structure for the untwisted Courant algebroid
$TM\oplus T^*M$. By Example \ref{Cf}, choosing $K=\cF\oplus N^*C$,
we know that we can reduce this Courant algebroid  and obtain the standard
Courant algebroid on $\un{C}$.

Using Prop. \ref{reddirac} now we show that $L:=D|_C$ also descends.  $L\cap \kp$ has 
constant rank since
it's isomorphic to $\cF^{\circ}$. To check \eqref{diracdes} we use the fact that $K$
is spanned by
  closed 1-forms  and   hamiltonian vector fields of functions
vanishing on $C$. The former  act trivially, the latter (acting by Lie derivative)
map $\Gamma(L)$ to itself because  hamiltonian vector fields
preserve the Poisson structure. An arbitrary section of $K$ maps  $\Gamma(L\cap \kp)]$
to  $\Gamma(L+K)$ by the ``Leibniz rule in the first entry'' (see Section \ref{review}), so 
 \eqref{diracdes} is satisfied.
   Further it's known \cite{Cou} that the integrability of $\Pi$ is equivalent to $\Gamma(D)$ being
closed under the Courant bracket, so \eqref{diracint} holds.
Hence Prop. \ref{reddirac}
tells us that $\Pi$ descends to a Dirac structure on $\un{C}$. This
of course is the Poisson structure  obtained by the classical coisotropic reduction.
\end{ex}

\section{On the reduced Courant algebroid}\label{Sev}
Using the methods of Section \ref{dirac}
we derive some results on the reduced Courant algebroid obtained in Thm. \ref{redcour}. 
In this section  $E$ is an exact Courant algebroid over $M$
 andlet $C$  a submanifold  endowed with a coisotropic
subbundle $\kp$ of $E$ satisfying $\pi(\kp)=TC$.

\subsection{Adapted splittings}
 
In this subsection we consider ``good'' splittings of an exact Courant algebroid $E\rightarrow M$,
and using their existence we determine simple data on a foliated submanifold  that induce
an exact Courant algebroid on the leaf space (Prop. \ref{wgs}).

\begin{defi}\label{adapted} We call a  splitting $\sigma:TM \rightarrow E$ of the
sequence $\eqref{exact}$ \emph{adapted to $K$} if

\begin{itemize}
\item[a)] The image of $\sigma$ is isotropic  
\item[b)]  $\sigma(TC)\subset \kp$
\item[c)] for any vector field $X$ on $C$ which is projectable to $\un{C}$
we have $\sigma (X)\in \Gbw$.
\end{itemize}
\end{defi}

\begin{remark}\label{K}
For such a splitting it follows automatically that
$\sigma(\pi(K))\subset K$. Indeed by $\pi(\kp)=TC$, b) in the
definition above and $\kp\cap T^*M=(\pi(K))^{\circ}$ we have
$\kp=\sigma(TC)+(\pi(K))^{\circ}$. Now $\la \sigma(\pi(K)),
\sigma(TC) \ra=0$ by a) in the definition above and
$\la \sigma(\pi(K)), (\pi(K))^{\circ} \ra=0$ 
Hence
$\sigma(\pi(K))$ has zero symmetric pairing with $\kp$.
\end{remark}

\begin{lemma}\label{speq}
By the prescription $\sigma \mapsto L_{\sigma}:=\sigma(TM)$,   
   splittings $\sigma$ adapted to $K$ correspond exactly to subbundles $L_{\sigma}\subset E$ with $\pi(L_{\sigma})=TM$ satisfying  
\begin{itemize}
\item[a)] $L_{\sigma}$  is maximal isotropic  
\item[b)] $\pi(L_{\sigma}\cap \kp)=TC$ 
\item[c)] $[\Gamma(K),\Gamma(L_{\sigma}\cap \kp)]\subset \Gamma(L_{\sigma}+K)$
 \end{itemize}
\end{lemma}
 \begin{proof} We just  show that c) is equivalent to item c) in 
Def. \ref{adapted}. For the implication ``$\Rightarrow$''
  we   have to check that if $k\in \Gamma(K)$ and $X$ is a projectable vector 
field on $C$ then $[k,\sigma(X)]\in \Gamma(K)$. Since $\sigma(X)\subset L_{\sigma}\cap \kp$
 this bracket is a section of $L_{\sigma}+K$. Further, since  $\pi([k,\sigma(X)])=[\pi(k),X]\subset \pi(K)$ 
 (because  $X$ is projectable) it actually lies in  $(L_{\sigma}+K)\cap \pi^{-1}(\pi(K))=K$.
 
The other implication  follows because $L_{\sigma}\cap \kp=\sigma(TC)$ admits a frame
of basic sections, namely $\sigma(X)$ as $X$ ranges over projectable vector fields on $C$.
\end{proof}

\begin{lemma}\label{newwgs}
Let $E$ be an exact Courant algebroid over a manifold $M$, $C$ a submanifold endowed
with a regular integrable foliation $\cF$ so that $C/\cF$ be smooth,
 and 
 $L$ a maximal isotropic subbundle $L\subset E|_C$ with
$\pi(L)=TC$ such that
$ [\Gamma(K),\Gamma(L)]\subset \Gamma(L)$
where $K:=L\cap \pi^{-1}(\cF)$.
Then  there exists a splitting adapted to $K$. \end{lemma}

\begin{proof}
Notice that $K$ is isotropic and has constant
rank, because $ker(\pi|_K)=K\cap T^*M=L\cap T^*M= N^*C$ has constant
rank and $\pi(K)=\cF$ has constant rank by assumption.  Also $\kp=L+\cF^{\circ}$, so $\pi(\kp)=TC$.
 Let
$\sigma: TM \rightarrow E$ be an isotropic splitting such that $\sigma(TC)\subset L$. Then $L_{\sigma}:=\sigma (TM)$ 
 clearly satisfies   conditions  a) and b) in Lemma \ref{speq}. Further it satisfies condition c)  because $L_{\sigma}\cap K^{\perp}+K= L$ (the inclusion ``$\subset$'' is easy, the other one follows from $K\cap T^*M=L\cap T^*M=N^*C$).
Hence, by Lemma   \ref{speq}, $\sigma$ is
 a splitting adapted to $K$. 
\end{proof}

The following proposition says that, with some regularity assumptions, splittings adapted to
$K$ exist if and only if the reduced exact Courant algebroid
$\un{E}$ as in Thm. \ref{redcour} exists. 
\begin{prop}\label{exists}
Let $K\rightarrow C$ be an isotropic subbundle of $E$ with $\pi(\kp)=TC$ and
assume that $\pi(K)$ be integrable and $\un{C}:=C/\pi(K)$ be smooth. Then  splittings adapted to $K$
exist if and only if  $\Gbw$ spans $\kp$ at every point of $C$.
\end{prop}
\begin{proof}
Assume first that a splitting $\sigma$ adapted to  $K$ exists.   Let $X$ be a projectable vector field on $C$. 
By c) of Def. \ref{adapted}  $\sigma(X)$ we will lie in $\Gbw$. Take a
function on $\un{C}$, pull it
back to a function  $C$ and  extend it to a
function $f$ on   $M$. Then $df|_C$ is a 
 section of $(\pi(K))^{\circ}=T^*M \cap \kp$. Further it lies in $\Gbw$: for any $k\in \Gamma(K)$
 we have
$$[k, df|_C]=-[ df|_C,k]+d\langle k, df|_C \rangle
\subset N^*C\subset K$$
 because $df$ as a closed 1-form acts trivially and it annihilates 
 $\pi(K)$.
Since $\kp=\sigma(TC)+(T^*M \cap \kp)$, taking all projectable vector fields $X$
and functions $f$ as above
  we see that
$\Gbw$ spans $\kp$ at every point of $M$.

Conversely, assume now that $\Gbw$ spans $\kp$ at every point of $MC$.
Then by Thm. \ref{redcour}
the Courant algebroid $\un{E}$ over
$\un{C}$ exists; let $\un{\sigma}:T\un{C} \rightarrow \un{E}$ be 
any isotropic splitting.
Denote by $L$ the preimage of the maximal isotropic subbundle $\un{\sigma}(T\un{C})$
under $p:\kp\rightarrow \kp/K \rightarrow \un{E}$. $L$ is a maximal isotropic subbundle
of $\kp$, and $\pi(L)=TC$. Furthermore, lifting sections of $\un{\sigma}(T\un{C})$ to basic sections of $\kp$, we see that  $\Gamma_{bas}(L)$ (the basic sections that lie in $L$)  spans $L$ at every point of $C$, hence from $[\Gamma(K), \Gamma_{bas}(L)]\subset \Gamma(K)$ we can conclude $[\Gamma(K), \Gamma(L)]\subset \Gamma(L)$.
Notice also that $L\cap \pi^{-1}(\pi(K))=p^{-1}(\un{\sigma}(T\un{C})\cap T\un{C})=K$.
Hence we can apply Lemma \ref{newwgs} and obtain a 
splitting of $E$ adapted to $K$. 
 \end{proof}

Putting together  
Lemma \ref{newwgs} , Prop. \ref{exists}  and Thm. \ref{redcour}
 we are obtain:
\begin{prop}\label{wgs}
Let $E$ be an exact Courant algebroid over a manifold $M$, $C$ a submanifold endowed
with a regular integrable foliation $\cF$ so that $C/\cF$ be smooth,
 and 
 $L$ a maximal isotropic subbundle $L\subset E|_C$ with
$\pi(L)=TC$ such that
$ [\Gamma(K),\Gamma(L)]\subset \Gamma(L)$
where $K:=L\cap \pi^{-1}(\cF)$.
Then  $E$ descends to an
exact Courant algebroid on $C/\cF$.
\end{prop}

\begin{proof}

  \end{proof}
 
\subsection{The \v{S}evera class of the reduced
Courant algebroid}
 
  In Theorem
\ref{redcour} we showed that, when certain assumptions are met, one
obtains an exact Courant algebroid $\un{E}$ over the quotient
$\un{C}$ of $C$ by the distribution $\pi(K)$. In this subsection we
will discuss how to obtain the \v{S}evera class of $\un{E}$ from the
one of $E$.

We start observing that  
  $j^*H_{\sigma}$
 descends to a 3-form on $\un{C}$,
where $j$ is the inclusion of $C$ in $M$.
Since $H_{\sigma}$ is closed
  we just need to check $i_X(j^*H_{\sigma})=0 $ where $X\in \pi(K)$. 
 Extend $X$ to a
vector field tangent to $\pi(K)$; take vectors $Y,Z \in T_pC$ and
extend them locally to projectable vector fields of $C$. Since
$\sigma$ is an splitting adapted to $K$
 we know that $\sigma(Y)\in \Gbw$, and since
 $\sigma(X)\subset K$ (by Remark \ref{K}) we have $[\sigma (X),\sigma
 (Y)]\subset K$.
  Therefore, indeed,
\begin{eqnarray} H_{\sigma}(X,Y,Z)&=& 2 \la [\sigma (X),\sigma (Y)], \sigma
(Z)\ra =   0.
\end{eqnarray}
Even more is true by the following, 
which is an analog of Prop. 3.6 of \cite{BCG} (but unlike that proposition
 does not involve equivariant cohomology; see also \cite{LT1,LT2}).

\begin{prop}\label{pavel}
Assume that $\un{C}$ is a smooth manifold. If $\sigma$ is  a
splitting adapted to $K$ then $j^*(H_{\sigma})$ descends to a closed 3-form on $\un{C}$ which
represents the \v{S}evera class of $\un{E}$.
\end{prop}
\begin{proof}
$L_{\sigma}:=\sigma(TM)$ satisfies the conditions listed in Lemma \ref{speq}; in particular condition b)   implies that $L_{\sigma}\cap \kp$ has constant rank, and condition c) is just eq. \eqref{diracdes}. Hence by Prop. \ref{reddirac} $L_{\sigma}$ descends to a maximal isotropic subbundle of $\un{E}$, which by condition b) is the image of a splitting $\un{\sigma}:T\un{C}\rightarrow \un{E}$.
If $X$ is a projectable vector field on $C$ and $\un{X}$ the corresponding vector field on $\un{C}$, then $\sigma(X)\in \Gamma(L_{\sigma}\cap \kp)$ projects to $\un{\sigma}(\un{X})$ under $\kp\rightarrow \un{E}$, as can be seen using the isomorphisms $L_{\sigma}\cap \kp\cong TC$ and $\un{\sigma}(T\un{C})$ given by the anchors.

To compute the 3-form on $\un{C}$ induced by $j^*H_{\sigma}$ pick
three tangent vectors on $\un{C}$ at some point $\un{p}$, which by
abuse by notation we denote by $\un{X},\un{Y},\un{Z}$. Extend them
to vector fields on $\un{C}$ and lift them   to obtain projectable vector fields
$X,Y,Z$. $\sigma(Z)$ lies in
$\Gamma_{bas}(\kp)$, and as seen above it is a lift of $\un{\sigma}(\un{Z})\in
\Gamma(\un{E})$. The same holds for $X$ and $Y$, therefore, by the
definition of Courant bracket   on
$\un{E}$, we know that $[\sigma(X),\sigma(Y)]\in \Gamma_{bas}(\kp)$
is a lift of $[\un{\sigma}(\un{X}),\un{\sigma}(\un{Y}) ] \in
\Gamma(\un{E})$. Hence
\begin{equation*} H_{\sigma}(X,Y,Z)= 2 \la [\sigma (X),\sigma (Y)], \sigma (Z)\ra
=  2 \la  [\un{\sigma}(\un{X}),\un{\sigma}(\un{Y}) ] ,
\un{\sigma}(\un{Z}) \ra.\end{equation*} That is, $H_{\sigma}$
descends to   the curvature 3-form of $\un{E}$ induced by the
isotropic splitting $\un{\sigma}$.
\end{proof}

\begin{remark}
If $\sigma$ and $\hat{\sigma}$ are any two isotropic splittings for
$E\rightarrow TM$ then there is a 2-form $b\in \Omega^2(M)$ for which
$\sigma(X)-\hat{\sigma}(X)=b(X,\cdot)\in T^*M$ for all $X\in TM$. It
is also known that $H_{\sigma}$ and $H_{\hat{\sigma}}$ differ by
$db$.
 Now let  $\sigma$ and $\hat{\sigma}$ be adapted to $K$. Then $j^*b$  descends to a 2-form on
 $\un{C}$.Indeed, if   $X\in\pi(K)$,
 $b(X,\cdot)=\sigma(X)-\hat{\sigma}(X)\in K\cap
 T^*M= N^*C$, so the interior product of $X$ with $j^*b$ vanishes, and the same is true for $d(j^*b)$ 
as the difference of 3-forms which descend to $\un{C}$. 
   This is consistent with the fact that by Prop. \ref{pavel}
$H_{\sigma}$ and $H_{\hat{\sigma}}$ descend to 3-forms that
represent the same element of  $H^3(\un{C},\RR)$  (namely the
\v{S}evera class of $\un{E}$).
\end{remark}

As an instance of how a splitting adapted to $K$ is used to compute
the \v{S}evera class of the reduced Courant algebroid
we   revisit Example 3.12 of \cite{BCG}.
\begin{ex}
Let $M=C=S^3\times S^1$, denote by $\partial_t$ the infinitesimal
generator of the action of the circle on $S^3$  giving rise to the
Hopf bundle $p:S^3\rightarrow S^2$, and by $s$ the coordinate on the
second factor $S^1$. Let $E=TM\oplus T^*M$ the untwisted (i.e.
$H=0$) Courant algebroid on $M$. We choose the rank-one  subbundle
$K$ to be spanned by $\partial_t+ds$. Choose a connection one form
$\alpha$ for the circle bundle $S^3\rightarrow S^2$, and denote by
$X^H\in TS^3$ the horizontal lift of a vector $X$ on $S^2$.
$K^{\perp}$ is spanned by $\{\partial_t,\partial_s-\alpha,X^H,p^*\xi,
ds\}$ where $X$ (resp. $\xi$) runs over all vectors (resp. covectors) on  $S^2$. Since
$ds$ is closed the adjoint action  of $\partial_t+ds$ is just the
Lie derivative w.r.t. $\partial_t$, which kills any of
$\partial_t,\alpha, X^H,p^*\xi, \partial_s,ds$. In particular
$\Gamma_{bas}(\kp)$ spans   $\kp$. Hence the assumptions of
Thm. \ref{redcour} are satisfied, and on $S^2\times S^1$ we have a
reduced exact Courant algebroid.
Now we choose the splitting $\sigma:TM\rightarrow K^{\perp}$ as
follows: 
$$\sigma(\partial_t)=\partial_t+ds,\;\;\;\sigma(X^H)=X^H+0 \text{ for
all }X\in TS^2,\;\;\; \sigma(\partial_s)=\partial_s-\alpha.$$ This
is a splitting adapted to $K$.
(to check that it maps projectable vector fields to elements of
$\Gamma_{bas}(\kp)$  use
$[\partial_t+ds,\cdot]=\cL_{\partial_t}$).

Now we compute $H_{\sigma}$. If $X,Y$
are vector fields on $S^2$ we have
$[\sigma(X^H),\sigma(Y^H)]=[X^H,Y^H]+0=([X,Y]^H-F(X,Y)\partial_t)+0$
where $F\in \Omega^2(S^2)$ is the curvature of $\alpha$. Also
$[\sigma(\partial_s),\sigma(X^H)]=0+p^*(i_{X}F)$, and the analog
computation for other other combinations of pairs of
$\sigma(\partial_t),\sigma(X^H),\sigma(\partial_s)$ is zero. From
this we deduce that  $H_{\sigma}=p^*F\wedge ds$, which descends to the 3-form
$F\wedge ds$ on $S^2\times S^1$, and by Prop. \ref{pavel} it
represents the \v{S}evera class of  
$\un{E}$.

As pointed out in \cite{BCG} $F\wedge ds$ defines a non-trivial
cohomology class. An ``explanation'' for this fact is that by
Prop. \ref{pavel} to obtain a 3-form on $C$ that descends to a
representative of the \v{S}evera class of $\un{E}$ we need to choose
a splitting adapted to $K$; the trivial splitting
$\hat{\sigma}$,  which delivers
$H_{\hat{\sigma}}=0$, fails to be one because it does not map into
$\kp$.
\end{ex}

\section{The case of generalized complex structures}\label{gcs}

Let $E$ be an exact Courant algebroid over $M$. Recall that a
generalized complex structure is a vector bundle endomorphism $\cJ$
of $E$ which preserves $\langle \cdot,\cdot \rangle$, squares to
$-Id_E$ and for which the Nijenhuis tensor
 \begin{equation}\label{nij}
N_{\cJ}(e_1, e_2) := [\cJ e_1,\cJ e_2] -[e_1, e_2] - \cJ([e_1,\cJ
e_2]+[\cJ e_1, e_2]).\end{equation} vanishes.

The analog of the following proposition when a group action is present
is Thm. 4.8 of \cite{SX}; we borrow the first part of our proof from them, 
but use different arguments to prove the integrability of the reduced generalized complex structure.

\begin{prop}[Generalized complex reduction]\label{redgcs}
Let $E\rightarrow M$ and  $K\rightarrow C$ satisfy the assumptions
of Thm. \ref{redcour}, so that we have an exact Courant algebroid
$\un{E}\rightarrow \un{C}$. Let $\cJ$ be a generalized complex
structure on $M$ such that $\cJ K \cap K^{\perp}$ has constant rank
and is contained in $K$. Assume further that 
$\cJ$ applied to any basic
section of $\cJ \kp \cap K^{\perp}$ is again a basic section. Then
$\cJ$ descends to a generalized complex structure $\un{\cJ}$ on
$\un{E}\rightarrow \un{C}$.
\end{prop}
\begin{remark}
The linear algebra conditions on  $\cJ K \cap K^{\perp}$ are in particular
satisfied when $\cJ(K)=K$, in which case $\cJ K \cap
K^{\perp}=K$. The opposite extreme  is when $\cJ K \cap
K^{\perp}=\{0\}$.
\end{remark}
\begin{proof}
First we show that $\un{J}$ induces a smooth\footnote{This is clear when $\cJ$ preserves $\kp$.} endomorphism of the
vector bundle
$\kp/K$ over $C$. Indeed $\cJ K\cap \kp\subset K$ is equivalent to
$\cJ \kp+K\supset \kp$, so that $\kp=\kp \cap (\cJ \kp +K)=(\kp\cap
\cJ\kp)+K$. From this it is clear that $\kpj$ maps surjectively
under $\Pi:\kp\rightarrow \kp/K$. Since $ker(\Pi|_{\kpj})= (\kpj)\cap
K=K\cap \cJ \kp$,  by our constant rank assumption we obtain a smooth
vector bundle
 $\kpj/ker(\Pi|_{\kpj})$ canonically isomorphic to  $\kp/K$.  

 We use again the
assumption $\cJ K \cap K^{\perp}\subset K$, interpreting it as
follows: if $e$ lies in the kernel of $\Pi:\kp\rightarrow \kp/K$ and 
$\cJ e \in \kp$ then
$\cJ e$ is still in the kernel. This applies in particular to all
$e\in ker(\Pi|_{\kpj})$ (since $\kpj$ is $\cJ$-invariant), so we
deduce that $\cJ$ leaves $ker(\Pi_{\kpj})$ invariant, i.e. $\cJ$
induces a well-defined endomorphism on $\kpj/ker(\Pi|_{\kpj})
\cong \kp/K$.  Further it is clear that it squares to $-1$ and
preserves the induced symmetric pairing on $\kp/K$.

 Now take
a section $\un{e}$ of $\un{E}$, lift it to a (automatically basic)
section $e$ of $\kpj$. Then by assumption  $\cJ e$ is again a basic section;
 this
shows that the endomorphism on $\kpj/ker(\Pi|_{\kpj})$ descends to an
endomorphism $\un{\cJ}$ of $\un{E}$.

To show that $\un{\cJ}$ is integrable one
can apply the 
proof of Thm. 6.1 in \ref{BCG}, where
 $\cJ$ is encoded in terms of complex Dirac structures.
 Alternatively one can show that the Nijenhuis tensor of $\cJ$ (which vanishes) is a lift of the Nijenhuis tensor of
 $\un{\cJ}$. The computation is straightforward except for the fact that $[e_1,\cJ e_2]+[\cJ e_1, e_2]$
 is a section of $\cJ \kp$ for all $e_1,e_2 \in \Gamma_{bas}(\kp)$, which follows using the Leibniz rule C4).
 \end{proof}

In Prop. \ref{redgcs} the condition 
that $\cJ$ preserve $\Gbj$
does not follow from the integrability of
$\cJ$ (see Ex. \ref{compl} below for an explicit example).
 In Section
\ref{branes} we will consider submanifolds $C$ for which the integrability of $\cJ$
does imply all the assumptions of Prop. \ref{redgcs}, in analogy to the case of  
coisotropic submanifolds in the Poisson setting.  

\begin{ex}[Complex foliations] \label{compl}
Take   $E$ to be the standard Courant algebroid
and $\cJ$ be given by a complex structure $J$ on $M$. Take $\cF$ to
be a real integrable distribution on $M$ preserved by $\cJ$ (so $J$
induces the structure of a  complex manifold on each leaf of
$\cF$) and $K=\cF\oplus 0$, so that $\un{M}:=M/\pi(K)=M/\cF$  be smooth. The
generalized complex structure $\cJ$ preserves $K$. If $\cJ$
mapped $\Gamma_{bas}(\kp)$ into itself\footnote{This is equivalent to saying that 
for any vector field $X$ on $M$ which is projectable the vector field 
$J(X)$ is also projectable.}
then by  Prop. \ref{redgcs} it would follow that $\un{M}$ would have
an induced generalized complex structure. Further, it would
necessarily correspond to an honest complex structure on $\un{M}$
that makes $M\rightarrow \un{M}$ into a holomorphic map.
However there are examples for which such a complex structure on
$\un{M}$ does not exist; in \cite{Wi} Winkelmann quotes an example
where $M$ is a twistor space of real dimension 6 and $\un{M}$ is the
4-dimensional torus.
 \end{ex}

\begin{ex} [Symplectic foliations]
Take  again $E$ to be the standard Courant algebroid
and $\cJ$ be given by a symplectic form  $\omega$ on $M$. Take $K=\cF$ to
be a real integrable distribution on $M$. One checks that $\cJ K\cap \kp$ is contained
in $K$ only if it is trivial, which is equivalent to saying that the leaves of $\cF$
are symplectic submanifolds. $\cJ$ maps basic sections of $\cJ \kp\cap \kp
=\cF^{\omega}\oplus \cF^{\circ}$ into basic
sections iff the hamiltonian vector field $X_{\pi^*f}$ is a projectable vector field
 for any function $f$, where $pr:M \rightarrow \un{M}:=M/\cF$. When this is the case
the induced generalized complex structure on $\un{M}$ is the symplectic structure
 given
by the isomorphism of vector spaces $\cF_x^{\omega}\cong T_{pr(x)}\un{M}$ (where
$x\in M$).
\end{ex}

\begin{remark}\label{nochar}  
It is known that a generalized complex manifold $(M,\cJ)$ comes with a canonical
Poisson structure $\Pi$, whose  sharp map $\sharp$ is given by the
composition $T^*M \hookrightarrow E \overset{\cJ}\rightarrow E
\overset{\pi} \rightarrow TM$. 
If in Prop. \ref{redgcs} we assume that $\cJ$ preserves
$K$, then $C$ is a necessarily a coisotropic submanifold, because
from $N^*C= (\pi(\kp))^{\circ}=K\cap
ker(\pi)\subset K$ we have $\sharp (N^*C)=\pi (\cJ N^*C)\subset \pi (K)\subset\pi (\kp)=TC$.
So $C/\sharp N^*C$ (if smooth) has an induced Poisson structure. 
We  know that also $\un{C}:=C/\pi(K)$ has a Poisson structure, induced from the 
 reduced generalized complex structure.
In general $\pi(K)$ is \emph{not} the characteristic distribution of $C$; we
just have an inclusion $\sharp N^*C \subset \pi(K)$\footnote{A case in
which this inclusion is strict is when $\cJ$ corresponds to the standard complex structure
on $M=\CC^n$ (with complex coordinates $z_k=x_k+iy_k$)
and $K=span\{\frac{\partial}{\partial x_1},\frac{\partial}{\partial y_1}\}$.}.
The Poisson bivector on $C/\pi(K)$ is just the pushforward of  the one on $M$.
\end{remark}

Given an exact Courant algebroid $E$ on $M$, recall that  a
\emph{generalized K\"ahler structure} consists of two commuting
generalized complex structures $\cJ_1,\cJ_2$  such that the
symmetric bilinear form on $E$ given by $\la \cJ_1\cJ_2\cdot,\cdot
\ra$ be positive definite. The following result borrows the proof of    Thm. 6.1 of
\cite{BCG}.
\begin{prop}[Generalized K\"ahler reduction]\label{redgk}
Let $E\rightarrow M$ and  $K\rightarrow C$ satisfy the assumptions
of Prop. \ref{redcour}, so that we have an exact Courant algebroid
$\un{E}\rightarrow \un{C}$. Let $\cJ_1,\cJ_2$ be a generalized
K\"ahler structure on $M$ such that $\cJ_1 K =K$. Assume further
that  $\cJ_1$ maps $\Gamma_{bas}(\kp)$ into itself and that $\cJ_2$
maps $\Gamma_{bas}(\cJ_2 \kp \cap K^{\perp})$ into itself.
 Then
$\cJ_1,\cJ_2$ descend to a generalized K\"ahler  structure on
$\un{E}\rightarrow \un{C}$.\end{prop}
\begin{proof}
By  Thm. \ref{redgcs} $\cJ_1$ induces a generalized complex structure
$\un{\cJ_1}$ on $\un{E}$. The orthogonal $K^{\cG}$ of $K$ w.r.t. $\la
\cJ_1\cJ_2\cdot,\cdot \ra$ is $(\cJ_2\cJ_1K)^{\perp}=\cJ_2 \kp$. Because of the identity
$\kp=K\oplus (K^{\cG}\cap \kp)$ the restriction to $\cJ_2 \kp \cap
\kp$ of the projection $\kp \rightarrow \kp/K$ is an isomorphism.
 So we can apply  Prop. \ref{redgcs} to $\cJ_2$  and obtain a
generalized complex structure $\un{\cJ_2}$ on $\un{E}$. Notice that
both $\cJ_1$ and $\cJ_2$ preserve $\cJ_2 \kp \cap \kp$; pulling back
sections of $\un{E}$ to basic sections of $\cJ_2 \kp \cap \kp$ one
sees that $\un{\cJ_1},\un{\cJ_2}$ form a generalized K\"ahler
structure on $\un{E}$.
\end{proof}
 
\section{The case of (weak) branes}\label{branes}

In this section we define branes and show that they admit a natural quotient
which is a generalized complex manifold endowed with a space-filling brane. Then we notice that quotients of more general objects, which we call ``weak branes'',
also inherit a generalized complex structure; examples of weak branes are coisotropic submanifolds in symplectic manifolds. Finally we show how weak branes 
can be obtained by passing from a generalized complex manifold
to a suitable submanifold.

\subsection{Reducing branes}

\begin{defi}\label{gs}
Let $E$ be an exact Courant algebroid over a manifold $M$. A
\emph{generalized submanifold} is a pair $(C,L)$ consisting of a
submanifold $C\subset M$ and a maximal isotropic subbundle $L\subset
E$ over $C$ with $\pi(L)=TC$ which is closed under the Courant
bracket
 (i.e.  $[\Gamma(L),\Gamma(L)]\subset \Gamma(L)$ with the conventions of Remark \ref{D}).
\end{defi}

This definition, which already appeared in the literature\footnote{It appeared in  Def. 3.2.2  of \cite{HU} with the name ``maximally isotropic extended
submanifold''. Also,
a subbundle $L$ as above but for which we just ask 
 $\pi(L)\subset TC$ is called
generalized Dirac structure in Def. 6.8 of \cite{AX}
(in the setting of the
skew-symmetric Courant bracket).},
 is just a splitting-independent rephrasing\footnote{Up to a sign,
since Def. 7.4 of \cite{Gu} requires $i^*H_{\sigma}=dF$ (in the notation of this lemma).} of
Gualtieri's original definition (Def. 7.4 of \cite{Gu}). See also Lemma 3.2.3 of \cite{HU}.	
\begin{lemma}\label{equiv}
Let $E$ be  an exact Courant algebroid over $M$. Choose an isotropic
splitting $\sigma$ for $E$, giving rise to an isomorphism of Courant
algebroids $(E,[\cdot,\cdot])\cong (TM\oplus
T^*M,[\cdot,\cdot]_{H_{\sigma}})$ where $H_{\sigma}$ is the
curvature 3-form of the splitting (see Section \ref{review}). Then
pairs $(C,L)$ as in Def. \ref{gs} correspond bijectively to pairs
$(C,F)$, where $F\in \Omega^2(C)$ satisfies
 $-i^*H_{\sigma}=dF$  (for $i$ the inclusion of $C$ in $M$).\\
\end{lemma}
\begin{proof}
The fact that $L\subset E$ is maximal isotropic and $\pi(L)=TC$
means that under the isomorphism it maps to
 $$\tau_C^F:=\{(X,\xi)\in TC\oplus T^*M|_C:
\xi|_{TC}=i_X F\}$$ for some 2-form $F$ on $C$. The correspondence
$L \leftrightarrow F$ is clearly bijective.
The integrability conditions correspond because of the following, which follows from a 
 straight-forward computation:   if $X_i+\xi_i$ are sections of $\tau_C^F$ then
\begin{equation}\label{SWsign}
2\langle [X_1+\xi_1,X_2+\xi_2], X_3+\xi_3 \rangle=(i^*H+dF)(X_1,X_2,X_3).
\end{equation}
\end{proof}

By Lemma \ref{equiv} the following definition is equivalent to
Gualtieri's original one (i.e. to Def. 7.6 of \cite{Gu}, again up to a sign):

\begin{defi}\label{gcsub}
Let $E$ be an exact Courant algebroid over a manifold $M$ and $\cJ$
be a generalized complex structure on $E$. A \emph{generalized
complex submanifold} or \emph{brane} is a
generalized submanifold $(C,L)$ satisfying $\cJ(L)=L$.
\end{defi}

Now we state the main theorem of this paper.  Recall that we gave
the definition of coisotropic submanifold in Remark \ref{nochar}.
 
\begin{thm}[Brane reduction]\label{redbrane}
Let $E$ be an exact Courant algebroid over a manifold $M$, $\cJ$ a
generalized complex structure on $E$, and $(C,L)$ a brane.
 Then $C$ is coisotropic w.r.t. the
Poisson structure induced by $\cJ$ on $M$. If the quotient $\un{C}$
of $C$ by its characteristic foliation
 is
smooth, then
\begin{itemize}
\item [a)]$E$ induces an exact Courant algebroid $\un{E}$ over
$\un{C}$
\item [b)] $\cJ$ induces a generalized complex structure $\un{\cJ}$ on $\un{E}\rightarrow \un{C}$
\item [c)] $L$ induces the structures of a space-filling brane on
$\un{C}$ and the \v{S}evera class of $\un{E}$ is trivial.
\end{itemize}
\end{thm}
 
\begin{proof}
Recall that the Poisson structure $\Pi$ induced by
$\cJ$ on $M$ (or rather its sharp map $\sharp$) is given by the
composition $T^*M \hookrightarrow E \overset{\cJ}\rightarrow E
\overset{\pi} \rightarrow TM$. Since $N^*C= (\pi(L))^{\circ}=L\cap
ker(\pi)\subset L$ we have $\sharp (N^*C)=\pi (\cJ  N^*C)\subset  \pi (L)=TC$, so $C$ is a coisotropic
submanifold. As above we let $\cF:=\sharp N^*C$, assume that it be a
regular distribution
 and that $\un{C}:=C/\cF$ be a smooth manifold.
 
\textbf{a)} $C$, $L$ and $\cF$ satisfy the assumptions  of Prop. \ref{wgs}. Hence    Thm. \ref{redcour} (with $K:=L\cap \pi^{-1}(\cF)$) delivers
  an exact Courant algebroid $\un{E}$ over
$\un{C}$.
  Notice that we have not made use of the integrability of $\cJ$ here, if not for the fact 
  that the induced bivector $\Pi$ is integrable and hence the distribution $\cF$  is involutive. 
 
\textbf{b)} Now we check that the assumptions of Prop. \ref{redgcs} are satisfied.
From $L\cap T^*M=N^*C$, the fact that $\cJ N^*C $ is contained in $L$
 and that it projects onto $\cF$  we deduce that  
$K=N^*C+\cJ N^*C$, which is clearly
  preserved by $\cJ$.
So we just need to check
that, for any basic section $e$ of $\kp$,
$\cJ e$ is again basic. 
Locally we can write $K=span\{(dg_i)|_C,\cJ (dg_i)|_C\}$ where $g_1,\dots,g_{codim(C)}$
are local functions on $M$ vanishing on $C$. Since each $dg_i$ is a closed one form,
$[(dg_i)|_C,\cJ e]\subset K$. Using the fact that the Nijenhuis tensor $N_{\cJ}$ vanishes \eqref{nij} we have
$$[\cJ (dg_i)|_C, \cJ e]=\cJ[\cJ (dg_i)|_C,  e]+\cJ[ (dg_i)|_C, \cJ e]
+[(dg_i)|_C,  e].$$
The first term on the r.h.s. lies in $K$ because $e$ is a basic section, and the last two because $dg_i$
is a closed 1-form. So $[\cJ (dg_i)|_C, \cJ e]\subset K$, hence $e$ is again a basic section.
Hence  the assumptions of   Prop. \ref{redgcs} are
satisfied, concluding the proof of b).

\textbf{c)}  We want to apply  Prop.
\ref{reddirac} to obtain a brane on $\un{C}$. Since $L\subset \kp$
the assumption \eqref{diracdes} needed for $L$ to descend reads
$[\Gamma(K),\Gamma(L)]\subset \Gamma(L)$, and the integrability
assumption \eqref{diracint} reads $[\Gamma_{bas}(L),\Gamma_{bas}(L)]\subset
\Gamma(L)$. As $L$ is closed under the bracket both assumptions
hold, and we obtain an (integrable) Dirac structure $\un{L}$ on
$\un{C}$. Furthermore from the fact that $\cJ$ preserves $L$ we see
that $\un{\cJ}$ preserves $\un{L}$. Hence $(\un{C},\un{L})$ is a brane for
the generalized complex structure $\un{\cJ}$ on $\un{E}$.

If we chose any isotropic splitting for $\un{E}$, as discussed in
Lemma \ref{equiv}, then $\un{L}$ gives rise to a 2-form
$\hat{F}$ on $\un{C}$ such that $-d\hat{F}$ equals the curvature of
the splitting, which hence is an exact 3-form. This concludes the proof of c)
and of the theorem.

\end{proof}

\begin{remark}\label{referee}
Let us denote with $(\un{C},\un{\cJ})$   a generalized complex manifold admitting a space-filling brane $(\un{C},\un{L})$. By Example 6.12 of \cite{Gu2}, using the splitting $T\un{C}\rightarrow \un{E}$ with image $\un{L}$ the 
exact Courant algebroid $\un{E}$ is identified with the untwisted Courant algebroid, and the
generalized complex structure assumes the form 
\begin{equation}\label{triangle}
\left(\begin{array}{cc} -\un{I}& \un{\Pi} \\0 & \un{I}^*\end{array}\right).
\end{equation}
Here $\un{I}$ is an 
honest complex structure on $\un{C}$, which corresponds to $\un{J}:\un{L}\rightarrow \un{L}$ under the anchor $\un{L}\cong T\un{C}$,
 and $\un{\pi}$ is the Poisson structure induced by $\un{\cJ}$.

Further  $\un{Pi}$ is the imaginary part a holomorphic Poisson bivectorfield on $\un{C}$.
Therefore Thm. \ref{redbrane} c) implies that the generalized complex structure $\un{\cJ}$ on the quotient of a brane (when smooth) is a holomorphic Poisson deformation of a complex structure.

Now assume the set-up of Thm. \ref{redbrane}, i.e. that $E$ is an exact Courant algebroid over   $M$, $\cJ$ a
generalized complex structure, and $(C,L)$ a brane with smooth quotient $\un{C}$. Encoding $\cJ$ by the  $+i$-eigenbundle of its complexification,
Gualtieri provides a direct way to describe  the reduced generalized complex structure $\un{J}$ we obtained in 
Thm. \ref{redbrane} (Cor. 6.6 of \cite{Gu2}).
 Let $\ell$ denote the $+i$-eigenbundle of $J:L\otimes \CC \rightarrow L\otimes \CC $,
and $A:=\pi(\ell)\subset TC\otimes \CC$. Then $A$ descends to a distribution $\un{A}$ on $\un{C}$ so that
$\un{A}\oplus \un{\bar{A}}= T\un{C}\otimes \CC$, which defines a complex structure on $\un{C}$. This agrees with the complex structure $\un{I}$ that appeared above, since
the $+i$-eigenbundle of the complexification of $\un{I}$ is the image under $\un{\pi}$ of the $+i$-eigenbundle of  $\un{\cJ}: \un{L}\otimes \CC \rightarrow \un{L}\otimes \CC$, which is just $\un{\pi(\ell)}=\un{A}$. Further the Poisson structure $\un{\pi}$ is obtained simply  by coisotropic reduction. Hence by means \eqref{triangle}we recover the reduced generalized complex structure $\un{\cJ}$. 
\end{remark}

\begin{remark}
We saw in Thm. \ref{redbrane} that branes $C$  are coisotropic and
their quotient by the characteristic foliation is endowed with a generalized complex
structure. As pointed out in Remark \ref{nochar}, if one starts with a
$\cJ$-invariant coisotropic
subbundle $\kp$ of $E|_C$ (instead of constructing one from the brane $(C,L)$ as in Thm.
\ref{redbrane}) in
general it is  a different quotient of $C$ that is endowed with a generalized  complex structure
(via Prop. \ref{redgcs}).
If one picks just any arbitrary coisotropic submanifold $C$, its quotient by the characteristic foliation 
inherits a Poisson structure, but in general  it does not inherit a generalized
 complex structure:
take for example any odd dimensional submanifold of a complex manifold. 
\end{remark}

\begin{remark}\label{holo}
When the characteristic foliation of a brane $(C,L)\subset M$ is regular, 
using coordinates adapted to the foliation one sees that
the quotient of small
enough open sets $U$ of $C$ by the characteristic foliation is smooth, 
and Thm. \ref{redbrane} gives a local statement.
However in general the characteristic foliation is singular, as the following example
shows. Take $M=\CC^2$, the untwisted exact Courant algebroid as $E$,
and as $\cJ$ take 
$\left(\begin{smallmatrix}   I & \Pi \\ 0 & -I^*
\end{smallmatrix} \right)$. Here $I(\partial_{x_i})=\partial_{y_i}$ is the canonical complex structure
on $\CC^2$ and $\Pi=y_1(\partial_{x_1}\wedge \partial_{x_2}-\partial_{y_1}\wedge
\partial_{y_2})-x_1
(\partial_{y_1}\wedge \partial_{x_2}+\partial_{x_1}\wedge
\partial_{y_2})$ is the imaginary part of the holomorphic Poisson bivector (see
\cite{Gu2}\cite{Hi})
$z_1\partial_{z_1}\wedge \partial_{z_2} $. It is easy to check that
$C=\{z_2=0\}$ with $F=0$ define a brane for $\cJ$, and that the characteristic
distribution of $C$ has rank zero at the origin and rank 2 elsewhere.
\end{remark}

\begin{ex}[Branes in symplectic manifolds] 
Consider a symplectic manifold $(M,\omega)$ and view it as a generalized complex structure $\cJ=\left(\begin{smallmatrix}  0 & -\omega^{-1} \\ \omega & 0
\end{smallmatrix}\right)$
on the standard Courant algebroid.
Let $(C,F)$ be a brane, i.e. $F$ is a closed 2-form on $C$ such that $L:=\tau_C^F$ is preserved by $\cJ$. $F$ descends to the quotient of $C$ by the characteristic distribution $\cF:=TC^{\omega}$, hence the isotropic subbundle $K=L\cap \pi^{-1}(\cF)$ defined in a) of the proof of Thm. \ref{redbrane} is just $\cF\oplus N^*C$. The reduced Courant algebroid $K^{\perp}/K$ is therefore canonically isomorphic to $T\un{C}\oplus T^*\un{C}$, and the reduced generalized  complex structure is $\un{\cJ}=\left(\begin{smallmatrix}  0 & -\un{\omega}^{-1} \\ \un{\omega} & 0
\end{smallmatrix}\right)$ where $\un{\omega}$ denotes the symplectic form on $\un{C}$ descending from  the coisotropic submanifold $C\subset (M,\omega)$.

The Dirac structure on $\un{C}$ obtained pushing forward $L$ is just $graph(\un{F})$, where $\un{F}\in \Omega_2(\un{C})$ denotes the pushforward of  $F\in \Omega_2(C)$.
Taking the tangent component of the  action of $\un{\cJ}$ on $graph(\un{F})$ delivers
$\un{I}:=- \un{\omega}^{-1}\un{F}$ for the induced complex structure on $\un{C}$ as in Remark \ref{referee}. Applying the gauge-transformation by the closed 2-form $\un{F}$ brings $\un{\cJ}$ into the form  \eqref{triangle} (with $\Pi=-\un{\omega}^{-1}$).

Further, as shown in 
Example 7.8 of \cite{Gu},  
$\un{F}+i\un{\omega}$
 is a holomorphic symplectic form on $\un{C}$.
\end{ex}

\begin{remark}
Suppose that in the setting of Thm.  \ref{redbrane} $E$  is
additionally endowed with some $\cJ_2$ so that $\cJ_1,\cJ_2$ form a
generalized K\"ahler structure. Then using  Prop. \ref{redgk} we see
that \emph{if} $\cJ_2$ descends to $\un{E}$  
then $\un{E}$ is endowed with a generalized K\"ahler structure too.
\end{remark}

\subsection{Reducing weak branes}

We weaken the conditions in the definition of brane; at least for the time being,
we refer to resulting object as ``weak branes''.
\begin{defi}\label{defweak}
Let $E$ be an exact Courant algebroid over a manifold $M$, $\cJ$ a
generalized complex structure on $E$.
We will call \emph{weak brane} a pair
 $(C,L)$  consisting of a submanifold $C$ and 
 a maximal isotropic subbundle $L\subset E|_C$ with
$\pi(L)=TC$ such that
\begin{eqnarray}\label{weak}
  \cJ(N^*C)\subset L ,\;\;\;\;\;\;\;\;\;\;\;\;\;\;
 [\Gamma(K),\Gamma(L)]\subset \Gamma(L)
\end{eqnarray}

(where $K:=L\cap \pi^{-1}(\cF)$ and $\cF:=\sharp N^*C$,
or equivalently $K=N^*C + \cJ N^*C$.)
\end{defi}

Notice  weak branes for which 
 $\cF$ has constant rank automatically satisfy the assumptions of Prop. \ref{wgs}.
 Also notice that in the proof of Thm. \ref{redbrane} (except for c)) we just
 used properties of weak branes,  hence we obtain
\begin{prop}\label{weakbrane}
If in Thm. \ref{redbrane} we let $(C,L)$ be a weak brane
then $C$ is a coisotropic submanifold and a) and b) of Thm. \ref{redbrane} still hold, i.e.
 there is a reduced Courant algebroid and a
reduced generalized complex structure on $\un{C}$ (when it is a smooth manifold).
\end{prop}

We describe how weak branes look like
in the split case, i.e. when
$E=(TM\oplus T^*M,[\cdot,\cdot]_H)$.
 We write 
 $\cJ$ 
in matrix form as $\left(\begin{smallmatrix}   A & \Pi \\ \omega & -A^*
\end{smallmatrix}\right)$ where $A$ is an endomorphism of $TM$, $\Pi$ the
Poisson bivector canonically associated to $\cJ$, and $\omega$ a
2-form on $M$.

\begin{cor}\label{description}
Let $C$ be a submanifold of $M$ and $F\in \Omega^2(C)$.
Fix  an extension   $B\in \Omega^2(M)$ 
of $F$. Then $(C,\tau_C^F)$ is a weak brane (with smooth
quotient $\un{C}$)
iff
$C$ is coisotropic (with smooth
quotient $\un{C}$), $A+\Pi B:TM\rightarrow TM$ preserves $TC$ , and the 3-form $dF+i^*H$ on $C$ descends to  $\un{C}$.

In this case the \v{S}evera class of the reduced Courant algebroid $\un{E}$ is represented by the
pushforward of
$dF+i^*H$. Further there is a splitting of $\un{E}$ in which the reduced generalized complex structure
is $$\tilde{\un{\cJ}}=\begin{pmatrix}  \tilde{\un{A}} & \tilde{\un{\Pi}} \\ \tilde{\un{\omega}} &
-\tilde{\un{A}}^*
\end{pmatrix},$$ where the endomorphism $\tilde{\un{A}}$ is the pushforward of $(A+\Pi B)|_{TC}$, 
the Poisson bivector $\tilde{\un{\Pi}}$ is induced by $\Pi$, and the 2-form $\tilde{\un{\omega}}$ is the pushforward
of $i^*(\omega -B\Pi B -BA - A^*B)$.
\end{cor}

\begin{proof}
Since   $K$ is $\tau_C^F \cap \pi^{-1}(\cF)$ 
eq. \eqref{SWsign} shows that $[\Gamma(K),\Gamma(\tau_C^F)]\subset \Gamma(\tau_C^F))$ is equivalent to the fact 
 that the closed 3-form 
 $i^*H+dF$ descend to $\un{C}$.
Now perform a $-B$-transformation; the
transformed objects are
 $\tilde{L}=TC\oplus N^*C$ and  
$\tilde{\cJ}=\left(\begin{smallmatrix}   \tilde{A} & \tilde{\Pi} \\ \tilde{\omega} & -\tilde{A}^*
\end{smallmatrix}\right)$,
with components $\tilde{A}=A+\Pi B$, $\tilde{\Pi}=\Pi$ 
and $\tilde{\omega}=\omega -B\Pi B -BA - A^*B$ (see for example \cite{Vared}). Hence we see that the first condition in \eqref{weak}
is equivalent to $C$ begin coisotropic and $A+\Pi B$ preserving $TC$ (a condition independent of the 
extension $B$). Further, since by the proof of Thm. \ref{redbrane} $\cJ$ preserves $TC\oplus \cF^{\circ}$
and $\cF\oplus N^*C$, it is clear that in the induced splitting of $\un{E}$ the components
of $\tilde{\un{\cJ}}$ are induced from those of $\tilde{\cJ}$.

Now we show that  the \v{S}evera class of the reduced Courant algebroid $\un{E}$ is represented by the
pushforward of $dF+i^*H$.  
By the proof of Prop.
 \ref{wgs} 
 any isotropic splitting $\sigma$ of $(TM\oplus T^*M,[\cdot,\cdot]_H)$ with $\sigma(TC)\subset \tau_C^F$ (for example one is given by
  $\sigma(X):=X+i_XB$)
 is automatically a splitting adapted to $K$. Hence by Prop. \ref{pavel}    
  $i^*H_{\sigma}$ pushes down to a representative of
 the \v{S}evera class of  $\un{E}$. Now   $i^*H_{\sigma}$ is just  $dF+i^*H$, 
 because for vectors $X_i\in TC$ we have
 $H_{\sigma}(X_1,X_2,X_3)=2\langle [\sigma(X_1),\sigma(X_2)],  \sigma(X_3)
 \rangle=(i^*H+dF)(X_1,X_2,X_3),$
 where we used $\sigma(X_i)\subset \tau_C^F$ and
     \eqref{SWsign} in the last equality.
\end{proof}

We use the characterization of Cor. \ref{description} in the following examples.
\begin{ex}[Coisotropic reduction]
If $\cJ$ corresponds to a symplectic structure on $M$, then any coisotropic submanifold $C$
endowed with $F=0$ is a weak brane. The generalized complex structure
on $\un{C}$ (assumed to be a smooth manifold) corresponds to the reduced symplectic form.
\end{ex}

If $\cJ$ corresponds to a complex structure, then any weak brane is necessarily a complex submanifold.
If $\cJ$ is obtained deforming a complex structure in direction of a holomorphic Poisson
structure \cite{Gu2}\cite{Hi} this is no longer the case, as in the following two examples.
In both cases however the reduced generalized complex structures we obtain
are quite trivial.

\begin{ex}
Similarly to Remark \ref{holo} take
  $M$ to be the open halfspace $\{(x_1,y_1,x_2,y_2):y_1> 0\}
  \subset
  \CC^2$, the untwisted exact Courant algebroid as $E$,
and as $\cJ$ take 
$\left(\begin{smallmatrix} I & \Pi \\ 0 & -I^*
\end{smallmatrix}\right)$ where $I(\partial_{x_i})=\partial_{y_i}$ is the canonical complex structure
on $\CC^2$ and $\Pi=y_1(\partial_{x_1}\wedge \partial_{x_2}-\partial_{y_1}\wedge
\partial_{y_2})-x_1
(\partial_{y_1}\wedge \partial_{x_2}+\partial_{x_1}\wedge
\partial_{y_2})$ is the imaginary part of the holomorphic Poisson bivector  
$z_1\partial_{z_1}\wedge \partial_{z_2}$. 
We now take $C=\{(x_1,y_1,x_2,0):y_1> 0\}$ and on $C$ the closed  2-form
$F:=-\frac{1}{y_1}dy_1\wedge dx_2$. We show that the pair $(C,F)$ forms a weak brane.
By dimension reasons $C$ is coisotropic (the characteristic distribution is regular and
spanned by $x_1\partial_{x_1}+
y_1\partial_{y_1}$), so we just have to check that $I+\Pi B$ preserves $TC$, where
$B$ the 2-form on $M$ given by the same formula as $F$. This is true as one computes
$I+\Pi B:\partial_{x_1}\mapsto \partial_{y_1},\;\;
\partial_{y_1}\mapsto -\frac{x_1}{y_1}\partial_{y_1},\;\;
\partial_{x_2}\mapsto -\frac{x_1}{y_1}\partial_{x_2}.$

Now we want to compute the generalized complex structure on $\un{C}$ given by Prop. \ref{weakbrane},
We do so by first applying the gauge transformation by $-B$ to obtain a generalized complex structure
$\tilde{\cJ}$ and then
using the diffeomorphism $\un{C}\cong (-\frac{\pi}{2},\frac{\pi}{2})\times \RR$
induced by $C\rightarrow (-\frac{\pi}{2},\frac{\pi}{2})\times \RR, (x_1,y_1,x_2)\mapsto
(\theta:=arctg(\frac{x_1}{y_1}),x_2)$.
The Poisson bracket of the coordinate functions $\theta$ and $x_2$ on $\un{C}$  is computed by
pulling back the two functions to $C$, extending them to the whole of $M$ and taking their Poisson bracket
there. This gives the constant function $1$.
Next the coordinate vector field $\partial_{\theta}$ on $\un{C}$ is lifted by the vector field
$\frac{x_1^2+y_1^2}{y_1}\partial_{x_1}$ on $C$, and of course $\partial_{x_2}$ on $\un{C}$
is lifted by $\partial_{x_2}$ on $C$. Applying the endomorphism $I+\Pi B$  of $TC$ 
we see the  induced endomorphism on $T\un{C}$ is just multiplication by $-tg(\theta)$.
Finally, the component $\tilde{\omega}$ of $\tilde{\cJ}$ is given by $-BI -B \Pi B -I^*B$,
which on $C$ restricts to the 2-form $\frac{1}{y_1^2}(y_1 dx_1-x_1 dy_1)\wedge dx_2$, which in turn is   
the pullback of the 2-form $(1+tg^2(\theta))d\theta\wedge dx_2$ on $\un{C}$. Hence  
 the induced generalized complex structure on $\un{C}$
 is $$\begin{pmatrix}  -tg(\theta)\cdot Id & \partial_{\theta}\wedge \partial_{x_2} \\ 
 (1+tg^2(\theta))d\theta\wedge dx_2 & tg(\theta)\cdot Id
\end{pmatrix}.$$
This is just the gauge transformation by the closed 2-form $tg(\theta)d\theta\wedge dx_2$ of the
generalized complex structure on $(-\frac{\pi}{2},\frac{\pi}{2})\times \RR$ that corresponds to
the symplectic form $d\theta\wedge dx_2$.
\end{ex}

\begin{ex} 
Similarly to the previous example we take 
  $M=  \CC^2$, the untwisted exact Courant algebroid as $E$,
and as $\cJ$ we take 
$\left(\begin{smallmatrix}  I & \Pi \\ 0 & -I^*
\end{smallmatrix}\right)$ where $I(\partial_{x_i})=\partial_{y_i}$ is the canonical complex structure
on $\CC^2$ and $\Pi=y_1(\partial_{x_1}\wedge \partial_{x_2}-\partial_{y_1}\wedge
\partial_{y_2})-x_1
(\partial_{y_1}\wedge \partial_{x_2}+\partial_{x_1}\wedge
\partial_{y_2})$.
Now we let $C$ be the   hypersurface  
$\{x_1^2+y_1^2=1\}$. The characteristic distribution is generated by $\partial_{y_2}$, so the quotient
$\un{C}$ is a cylinder.
Let $a,b,c \in C^{\infty}(C)$ so that, denoting by $F_{(a,b,c)}$ the pullback to $C$ of $$B_{(a,b,c)}:=
a\cdot dx_1 \wedge dy_1+b\cdot dx_1 \wedge dx_2+
c\cdot dy_1 \wedge dx_2 -y_1\cdot dx_1 \wedge dy_2+
+x_1\cdot dy_1 \wedge dy_2,$$
$dF_{(a,b,c)}$ descends\footnote{This happens exactly when $F_{(a,b,c)}$ is closed.}
  to $\un{C}$. One checks that $I^*+ B_{(a,b,c)}\Pi$ preserves $N^*C$, so that $(C,F_{(a,b,c)})$ is a weak
  brane. A computation analog to the one of the previous example shows that
  the reduced generalized  complex structure on $\un{C}=S^1\times \RR$ with coordinates $\theta$ and
  $x_2$ is given by 
   $$\begin{pmatrix}  \lambda_{(a,b)}\cdot Id & \partial_{\theta}\wedge \partial_{x_2} \\ 
 (1+\lambda_{(a,b)}^2)d\theta\wedge dx_2 & -\lambda_{(a,b)}\cdot Id
\end{pmatrix}$$
 where $\lambda_{(a,b)}\in C^{\infty}(\un{C})$ is the function that lifts to $-by_1+cx_1\in 
 C^{\infty}(C)$ via $C\rightarrow \un{C}$.
 Again this is a gauge transformation of the standard symplectic 
 structure on $S^1\times \RR$. 
 
 A consequence is that
 for no choice of $a,b,c$ as above the weak brane $(C,F_{(a,b,c)})$ is actually a brane.
Indeed if this was the case by Thm. \ref{redbrane} we would obtain a space-filling brane for a symplectic structure
 on $S^1\times \RR$; applying again Thm. \ref{redbrane}, by  Example 7.8 of \cite{Gu}, we would
 obtain the structure of a holomorphic symplectic manifold on $S^1\times \RR$, which can not exist
 because  
 holomorphic symplectic manifolds have real dimension $4k$ for some integer $k$.
\end{ex}

\subsection{Cosymplectic submanifolds}
Recall that a submanifold $\tilde{M}$ of a Poisson manifold $(M,\Pi)$ is 
	\emph{cosymplectic} if $\sharp N^*\tilde{M} \oplus T \tilde{M}=TM|_{\tilde{M}}$.
It is known (see for example \cite{Xudir}) that a cosymplectic submanifold inherits canonically a Poisson structure. 
The following lemma, which follows also from more general results of \cite{BS}, says that generalized complex structures are also inherited by cosymplectic submanifolds:
\begin{lemma}\label{cosympl}
Let $E$ be an exact Courant algebroid over a manifold $M$, $\cJ$ a
generalized complex structure on $E$ and $\tilde{M}$ 
a cosymplectic submanifold of $M$ (w.r.t. the natural Poisson structure
on $M$ induced by $\cJ$).
 Then $\tilde{M}$ is naturally endowed with a generalized complex
structure.
\end{lemma} 
\begin{proof}
We want apply Prop. \ref{redgcs} with $K=N^*\tilde{M}$ (so $\kp=\pi^{-1}(T \tilde{M})$).
 The intersection
$\cJ K\cap \kp$ is trivial. Indeed if $\xi \in N^*\tilde{M}$ and $\pi (\cJ \xi)\in  T\tilde{M}$ then by the definition of cosymplectic submanifold
 $\pi (\cJ \xi)=0$ (recall that $\sharp = \pi \cJ|_{T^*M}$)
and the restriction $\sharp$
 to $N^*\tilde{M}$ is injective, so that $\xi=0$. Further  all sections of $\kp$ are basic, so $\cJ$ maps the set of basic sections of $\cJ \kp\cap \kp$ into 
 itself. Hence the assumptions of Prop. \ref{redgcs} are satisfied and we obtain a 
generalized complex
structure on $\tilde{M}$. 
\end{proof}

Now we describe how a pair $(C,L)$ which doesn't quite satisfy the conditions 
of Def. \ref{defweak}
  can be regarded as a weak brane by passing to a cosymplectic submanifold.

\begin{prop}\label{prepois}
Let $E$ be an exact Courant algebroid over a manifold $M$, $\cJ$ a
generalized complex structure on 
$E$, $C$ 
a submanifold and $L$ a maximal isotropic  subbundle of $E|_C$ with $\pi(L)=TC$. Suppose that
 $\cJ(N^*C)\cap \pi^{-1}(TC)$ is contained in $L$ and has constant rank.
 Then there exists a submanifold $\tilde{M}$ (containing $C$) 
which inherits a generalized complex structure $\tilde{\cJ}$
 from $M$, and so that $\tilde{L}$ satisfies
$\tilde{\cJ}(\tilde{N}^*C)\subset \tilde{L}$. Here $\tilde{L}$ is the pullback of $L$ to $\tilde{M}$
and $\tilde{N}^*C$ the conormal bundle of $C$ in $\tilde{M}$.

Further assume that $[\Gamma({L}\cap {\pi}^{-1}(\cF)),\Gamma({L})]\subset \Gamma({L})$ where 
 $\cF:=\sharp N^*C\cap TC$ is the characteristic distribution of $C$.
Then
$ [\Gamma(\tilde{L}\cap \tilde{\pi}^{-1}(\cF)),\Gamma(\tilde{L})]\subset \Gamma(\tilde{L})$. 
Hence $(C,\tilde{L})$ is a weak brane in $(\tilde{M}, \tilde{\cJ})$.
 \end{prop}

 \begin{proof}
 Since the intersection of
 $\cJ(N^*C)$ and $\pi^{-1}(TC)$ has constant rank the same holds for their sum and for  $\pi 
 (\cJ(N^*C)+ \pi^{-1}(TC))=\sharp N^*C+TC$. Hence $C$ is a \emph{pre-Poisson} submanifold \cite{CZ}
 of $(M,\Pi)$. Fix  any 
complement $R$ of $\sharp N^*C+TC$ in $TM|_C$; by Theorem 3.3 of \cite{CZ}, ``extending'' $C$ in direction of
 $R$ we obtain a  submanifold $\tilde{M}$ of $M$ which is cosymplectic. By Lemma \ref{cosympl} we know
 that $\tilde{M}$ is endowed with a generalized complex
structure $\tilde{\cJ}$. 
Further by the same lemma $\cJ K\cap \kp$ is trivial.
The projection  $\kp \rightarrow \kp/K$ (for $K=N^*\tilde{M}$) maps $\cJ \kp\cap \kp$
isomorphically onto $\kp/K$, and $\tilde{\cJ}$ is induced by the action of $\cJ$ on 
$\cJ \kp\cap \kp$. Therefore, denoting by    $\tilde{L}:=L/K$ the pullback of $L$ to $\tilde{M}$,
 requiring 
$\tilde{\cJ}(\tilde{N}^*C)\subset \tilde{L}$ is equivalent to requiring that $\cJ(N^*C \cap (\cJ \kp \cap \kp))$
maps into $\tilde{L}$ under $\kp \rightarrow \kp/K$, which in turn means 
$\cJ(N^*C) \cap \kp
\subset L$. Now using $\kp=\pi^{-1}(T \tilde{M})$,
$T\tilde{M}|_C=R\oplus TC$ and recalling that $R$ was chosen so that
$R\oplus (\sharp N^*C + TC)=TM|_C$, it follows that 
$\cJ(N^*C) \cap \kp=\cJ(N^*C) \cap \pi^{-1}(TC)$. So our assumption ensures that
$\tilde{\cJ}(\tilde{N}^*C)\subset \tilde{L}$.

Finally notice that the projection $\kp \rightarrow \kp/K$ maps $L$ onto $\tilde{L}$.  Since
 $\pi^{-1}(\cF)$ is mapped onto  $\tilde{\pi}^{-1}(\cF)$  we also have that $L\cap \pi^{-1}(\cF)$  is mapped onto 
 $\tilde{L}\cap \tilde{\pi}^{-1}(\cF)$. Hence our assumption
  $[\Gamma(L\cap \pi^{-1}(\cF)),\Gamma({L})]\subset \Gamma({L})$ 
  implies
  $ [\Gamma(\tilde{L}\cap \tilde{\pi}^{-1}(\cF)),\Gamma(\tilde{L})]\subset \Gamma(\tilde{L})$.
  \end{proof}


\bibliographystyle{habbrv}
\bibliography{bibbrane}
\end{document}